\documentclass{amsart}
\usepackage[T1]{fontenc}
\usepackage{lmodern,microtype}
\usepackage{amsmath,amssymb,amsthm,mathtools}
\usepackage{xcolor}
\usepackage{tikz}
\definecolor{linkblue}{RGB}{20,45,85}
\usepackage[colorlinks=true,linkcolor=linkblue,citecolor=linkblue,urlcolor=linkblue]{hyperref}

\allowdisplaybreaks[1]
\numberwithin{equation}{section}
\newtheorem{theorem}{Theorem}[section]
\newtheorem{lemma}[theorem]{Lemma}
\newtheorem{corollary}[theorem]{Corollary}
\newtheorem{proposition}[theorem]{Proposition}
\theoremstyle{remark}
\newtheorem{remark}[theorem]{Remark}
\newcommand{\Pp}{\mathcal P}
\newcommand{\T}{\mathcal T}
\newcommand{\R}{\mathbb R}
\newcommand{\norm}[1]{\left\lVert#1\right\rVert}
\newcommand{\dual}[2]{\norm{#1}_{#2'}}
\newcommand{\diver}{\operatorname{div}}

\newcommand{\osc}{\operatorname{osc}}
\newcommand{\tr}{\operatorname{tr}}
\newcommand{\RT}{\mathcal{RT}}
\title[$p$-robust saturation on triangular meshes]{Polynomial-degree-robust saturation on triangular meshes}

\author[Kexin Ding]{Kexin Ding}
\address{School of Mathematical Sciences, Zhejiang University, 866 Yuhangtang Road, Hangzhou, Zhejiang 310058, People's Republic of China.}
\email{dingkexin@zju.edu.cn}

\author{Yuwen Li}
\address{School of Mathematical Sciences, Zhejiang University,
866 Yuhangtang Road, Hangzhou 310058, Zhejiang, People's Republic of China}
\email{liyuwen@zju.edu.cn}
\date{}
\begin{document}

\begin{abstract}
We prove a polynomial-degree-robust ($p$-robust) saturation property for conforming finite element approximations of the Poisson equation on shape-regular triangular meshes. Previously, a $p$-robust saturation estimate for the spectral method on a square was proved in [C. Canuto et al., A saturation property for the spectral-Galerkin approximation of a Dirichlet problem in a square, ESAIM Math. Model. Numer. Anal., 53 (2019), pp. 987--1003]. For every fixed $\lambda>1$, increasing the polynomial degree from $p$ to any $q\geq\lceil\lambda p\rceil+3$ reduces the energy error by a factor $\alpha_\lambda\in(0,1)$ independent of $p$ and the mesh size, up to a data oscillation term.
\end{abstract}
\maketitle

\section{Introduction}\label{sec:statement}
Let $\Omega\subset\R^2$ be a bounded polygonal Lipschitz domain and let
$\T$ be a conforming shape regular
triangulation of $\Omega$.
We consider the Poisson problem
\begin{equation*}
 -\Delta u=f\quad\text{in }\Omega,\qquad u=0\quad\text{on }\partial\Omega,
\end{equation*}
in which $f\in L^2(\Omega)$. We write $\Pp_r(T)$ for polynomials of total degree at most $r$ on $T$,
set $\Pp_{-1}=\{0\}$, and define
\begin{equation}\label{eq:spaces}
 V_r=\{v\in H^1_0(\Omega):v|_{T}\in\Pp_r(T)\text{ for all }T\in\T\},
 \qquad r\ge1.
\end{equation}
Let $(\cdot,\cdot)_\Omega$ and $\norm{\cdot}_\Omega$ denote the
$L^2(\Omega)$ inner product and norm, respectively, and set
$a(v,w)=(\nabla v,\nabla w)_\Omega$.
The exact and Galerkin solutions $u\in H^1_0(\Omega)$ and
$u_r\in V_r$~satisfy
\begin{equation}\label{eq:galerkin}
 \begin{aligned}
 a(u,v)&=(f,v)_\Omega,&& v\in H^1_0(\Omega),\\
 a(u_r,v_r)&=(f,v_r)_\Omega,&& v_r\in V_r.
 \end{aligned}
\end{equation}

Let $p$ denote the base polynomial degree and $q>p$ be an enriched polynomial degree. The energy error is
$E_r=\norm{\nabla(u-u_r)}_\Omega$. For each fixed $\lambda>1$, each integer $p\ge1$, and every
$f$ with $f|_T\in\Pp_{p-1}(T)$ on all $T\in\T$, we prove the following $p$-robust saturation under proportional degree enrichment 
\begin{equation}\label{eq:objective}
 E_q\le\alpha_\lambda E_p,\qquad \alpha_\lambda<1,\qquad
 \text{for }q\ge\lceil\lambda p\rceil+3,
\end{equation}
where $\alpha_\lambda\in(0,1)$ depends only on $\lambda$,
shape regularity, and $\Omega$. For general $f\in L^2(\Omega)$,
Theorem~\ref{thm:osc} will add an oscillation term to \eqref{eq:objective}.

The saturation assumption originates from a posteriori analysis in the
Bank--Weiser local Neumann estimators, whose reliability relies on
energy-error reduction in a richer space \cite{BW85}. 
Saturation also underlies the reliability analysis of hierarchical basis
estimators \cite{B96} and smoother-type
estimators \cite{LS26}.
Saturation itself holds for $P_1$--$P_2$ enrichment under sufficiently
small relative data oscillation \cite{DN02}, and for
$h$--$h/2$ refinement of triangular meshes up to patchwise data oscillation, subject to a computable test
\cite{CGG16}. On the other hand, interpolation versions of saturation property have been studied in \cite{BY15,BO17,BXY26}.

When $\Omega$ is a square, let $w_r$ be its Legendre spectral-Galerkin approximation for the Poisson problem in
$H^1_0(\Omega)\cap\Pp_r(\Omega)$. For each $\mu>1$ and $f\in\Pp_p(\Omega)$, Canuto,
Nochetto, Stevenson, and Verani \cite[Theorem~1.1]{CNSV19} prove the saturation estimate
\begin{equation*}
 m>\max\{\mu p,p+4\}
 \quad\Longrightarrow\quad
 \norm{\nabla(u-w_m)}_\Omega\le\beta_\mu\norm{\nabla(u-w_p)}_\Omega,
\end{equation*}
where $\beta_\mu\in(0,1)$ is independent of $p$, $m$, and $f$. 
Their experiments suggest the $p$-uniform
saturation would fail if the degree increment $m-p$ was fixed. The result in \cite{CNSV19} was generalized to  rectangular meshes with hanging nodes  ~\cite{W20}.

For triangular meshes, Bank, Parsania, and Sauter~\cite[Theorem~6.2]{BPS13} established a saturation estimate for $hp$-FEM
\begin{equation*}
 \norm{\nabla(u-u_{p+1})}_\Omega
 \le \kappa_p\norm{\nabla(u-u_p)}_\Omega,
\end{equation*}
where the contraction factor $\kappa_p<1$ approaches one as $p$ increases.
For adaptive $h$-refinement at fixed degree, an a posteriori
criterion in \cite{CDGV26} guarantees $p$-uniform error contraction.
For $p$-enrichment on triangular meshes, Canuto
et al.~\cite{CNSV17} showed that $p$-robust energy 
error contraction could be reduced to local saturation, that is, three dual norm comparisons on a reference
triangle~\cite[Theorem~6.3]{CNSV17}. Their experiments support local saturation, while a theoretical validation on triangular meshes remained open.

We prove the three dual norm comparisons on the reference triangle by
constructing polynomial operators that preserve the required edge and
volume moments. For each fixed enrichment factor $\lambda>1$, the
stability bounds for these operators are uniform in $p$. The three comparisons yield the vertex-patch residual estimate in
Proposition~\ref{prop:patch}, a key step towards a proof of $p$-robust saturation. The analysis uses two established polynomial lifting results. First, on a triangle,
the degree-preserving trace lifting from $H^{1/2}(\partial T)$ to $H^1(T)$ in \cite{AD}
is bounded independently of the polynomial degree. Second, the regularized
Poincar\'e operators in  \cite{CM} yield a $p$-robust right
inverse of divergence.

The next theorem is our main result.
\begin{theorem}\label{thm:main}
Let $f|_{T}\in\Pp_{p-1}(T)$ for each $T\in\T$.
For every $\lambda>1$, there is $\alpha_\lambda\in(0,1)$, depending only
on $\lambda$, shape regularity, and the domain, such that for every
integer $p\ge1$,
\begin{equation*}
q\ge\lceil\lambda p\rceil+3
 \quad\Longrightarrow\quad E_q\le\alpha_\lambda E_p.
\end{equation*}
\end{theorem}


For general data, let $\pi^{T}_{p-1}$ be the $L^2(T)$ projection onto
$\Pp_{p-1}(T)$ and set
\begin{equation}\label{eq:osc}
 \osc_p(f,\T)^2=
 \sum_{T\in\T}\left(\frac{h_{T}}{p}\right)^2
 \norm{f-\pi^{T}_{p-1}f}_{T}^2,\qquad h_{T}=\operatorname{diam}T.
\end{equation}

\begin{theorem}\label{thm:osc}
For every fixed $\lambda>1$, $f\in L^2(\Omega)$, $p\ge1$, and
$q\ge\lceil\lambda p\rceil+3$, we have
\begin{equation*}
 E_q\le\alpha_\lambda E_p+(1+\alpha_\lambda)C_{\rm osc}\osc_p(f,\T),
\end{equation*}
where $\alpha_\lambda$ is as in Theorem~\ref{thm:main} and $C_{\rm osc}$ is
independent of $\lambda$, $p$, and the mesh sizes, and depends only on shape
regularity.
\end{theorem}

Theorem~\ref{thm:osc} follows from Theorem~\ref{thm:main} by projecting
the data onto $\Pp_{p-1}$ and applying the local approximation estimate
\begin{equation}\label{eq:local-projection}
 \norm{v-\pi^{T}_{p-1}v}_{T}\le C\frac{h_{T}}{p}\norm{\nabla v}_{T},
 \qquad v\in H^1(T).
\end{equation}

In addition, the associated patch saturation estimate in Proposition~\ref{prop:patch} verifies
the local hypothesis of~\cite[Proposition~4.1]{CNSV17}. For resolved
polynomial data, one step of proportional $p$-enrichment on
D\"orfler-marked vertex patches therefore reduces the global energy
error by a factor strictly below one, uniformly in the mesh sizes and
polynomial degrees. Hence the \textsc{Reduce} module of $hp$-AFEM in~\cite{CNSV17optimal} can
achieve any fixed error reduction without coarsening, in a number of
enrichment steps independent of $p$. 

\begin{remark}[A posteriori error estimates and saturation]
The classical proof of $\mathcal{P}_1$--$\mathcal{P}_2$ saturation \cite{DN02} makes use of a residual-type error estimator as a bridge.
For $p$-robust saturation analysis, it is natural to use a $p$-robust equilibrated flux
estimator \cite{BPS09,EV15} intead. Solving mixed problems on vertex patches
gives a guaranteed upper bound for $E_p$ with a local efficiency constant
independent of $p$.

The additional step is
to bound $E_p$ by $\norm{\nabla(u_q-u_p)}_\Omega$ uniformly in $p$. The latter norm is the dual norm of the
residual of $u_p$ restricted to $V_q$. To prove the bound, we combine
successive edge corrections around each vertex, as in the equilibration
argument in~\cite{BPS09,EV15,DL26}, with stable polynomial
tests that preserve the volume and edge moments of the residual, see
Proposition~\ref{prop:patch}.
\end{remark}

The rest of the paper is organized as follows.
Sections~\ref{sec:boundary} and~\ref{sec:bubble} construct the boundary
and bubble moment operators. Subsection~\ref{sec:comparison} compares
the bubble argument with the square-domain proof of \cite{CNSV19}.
Section~\ref{sec:mixed} establishes the three reference estimates, and
Section~\ref{sec:patch} combines them with a polynomial divergence
lifting to obtain patch norming. Multiplication by a
vertex hat function raises the test degree by one. Section~\ref{sec:global}
proves saturation and the data-oscillation bound.

\section{Boundary moment preservation}\label{sec:boundary}

Throughout the paper, $C$ denotes a generic positive constant that may change from line to line. Unless stated otherwise, it is independent of the
polynomial degrees, mesh sizes, and enrichment parameters, but may depend on the fixed domain and mesh shape regularity.

The space $H^1_0(\Omega)$ is equipped with the energy norm $\norm{\nabla v}_\Omega$, which inducing the dual norm on $H^1_0(\Omega)'$.

\subsection{Local norms and moment preservation}\label{sec:prelim}
For a triangle $T$ and a nonempty union $D$ of its edges, we use
\begin{equation*}
 H_D^1(T)=\{v\in H^1(T):\tr_Dv=0\},\qquad
 \dual{F}{H_D^1(T)}=\sup_{0\ne v\in H_D^1(T)}
 \frac{|F(v)|}{\norm{\nabla v}_T}.
\end{equation*}
The homogeneous Dirichlet condition $\tr_Dv=0$ is the essential
boundary condition on $D$. We call the edges in $D$ essential edges.
When $D$ is empty, we use $H^1(T)/\R$, represented by functions of zero
volume mean. A functional on this quotient must annihilate constants.
Discrete dual norms use the same gradient norm restricted to polynomial
test functions. A volume function $F$ is identified with the functional
$v\mapsto(F,v)_T$.
For linear subspaces $V\subset U\subset X$, where $X$ has inner product
$(\cdot,\cdot)_X$, we write
$U\ominus_X V:=\{u\in U:(u,v)_X=0\ \forall v\in V\}$.

\begin{lemma}\label{lem:fortin}
Let $X$ be a normed space and $X_r\subset X$ a polynomial subspace.
Suppose that $J:X\to X_r$ is linear,
$\norm{Jv}_X\le C\norm v_X$, and $F(Jv)=F(v)$ for every functional
$F$ in a prescribed family. Then
\begin{equation*}
 \dual F X\le C\dual F{X_r}
\end{equation*}
for every functional in that family.
\end{lemma}
\begin{proof}
For $v\in X$, we have
$|F(v)|=|F(Jv)|\le\dual F{X_r}\norm{Jv}_X
\le C\dual F{X_r}\norm v_X$. Taking the supremum completes the proof.
\end{proof}

For $K=(0,1)$, we equip $H^{1/2}(K)$ with its Sobolev--Slobodeckij norm and use
\begin{equation*}
 \norm g_{H^{1/2}_{00}(K)}^2\simeq
 \norm g_{L^2(K)}^2+
 \int_K\int_K\frac{|g(s)-g(t)|^2}{|s-t|^2}\,ds\,dt+
 \int_0^1\frac{|g(s)|^2}{s(1-s)}\,ds.
\end{equation*}
If zero extension is required only across the endpoint $s=0$, respectively
$s=1$, the last weight is $1/s$, respectively $1/(1-s)$.
For an open chain of two edges, we use the same spaces in arclength
coordinates. Zero extension identifies these spaces with traces that
vanish on the complementary edges. We use the Sobolev trace theorem and
the interpolation identities
\begin{equation}\label{eq:trace-interpolation}
 [L^2(K),H^1_0(K)]_{1/2}=H^{1/2}_{00}(K),\qquad
 [L^2(\mathbb S^1),H^1(\mathbb S^1)]_{1/2}=H^{1/2}(\mathbb S^1).
\end{equation}
Here $\mathbb S^1=\R/\mathbb Z$ denotes the one-dimensional torus, and
its Sobolev spaces are understood periodically.

\subsection{The weighted Legendre filter}
On $I=(-1,1)$ let $L_j$ be the Legendre polynomial normalized by $L_j(1)=1$,
so $\norm{L_j}_I^2=2/(2j+1)$.
Let $P_m$ be the $L^2$ projector onto
$\Pp_m(I)$. Given $n\ge1$ and $r\ge n+2$, put
\begin{equation*}
 k=n+\left\lfloor\frac{r-n}{2}\right\rfloor-1,\qquad
 Z=\sum_{m=n}^{k}(m+1)
 =\frac{(k-n+1)(n+k+2)}2.
\end{equation*}
In particular $k\ge n$ and $r-k\ge2$. Define
\begin{equation*}
 S_{n,k}=\frac1Z\sum_{m=n}^{k}(m+1)P_m.
\end{equation*}
It maps into $\Pp_k$, is an $L^2$ contraction, and preserves every
$\Pp_n$ moment. Let $\omega_m=(m+1)/Z$ for $n\le m\le k$ and zero
otherwise. If $v=\sum_{j\ge0}\widehat v_jL_j$, then
\begin{equation*}
 S_{n,k}v=\sum_{j\ge0}s_j\widehat v_jL_j,\qquad
 s_j=\sum_{m\ge j}\omega_m.
\end{equation*}

We first apply the weighted average $S_{n,k}$ to preserve the
required moments and to control the derivative norm and endpoint errors. Its
output has degree at most $k$, so we use the remaining degrees
$k+1,\ldots,r$ to correct the endpoint values.

\begin{lemma}\label{lem:sn}
For $v\in H^1(I)$,
\begin{align*}
 \norm{(S_{n,k}v)'}_I
 &\le\left(1+\sqrt{\frac{2k(k+1)}Z}\right)\norm{v'}_I,\\
 |(v-S_{n,k}v)(\pm1)|
 &\le\sqrt{\frac2Z}\norm{v'}_I.
\end{align*}
\end{lemma}
\begin{proof}
\begingroup
For polynomial $v=\sum_jc_jL_j$ with $v'=\sum_jb_jL_j$, we have
\begin{equation*}
 c_j=\frac{b_{j-1}}{2j-1}-\frac{b_{j+1}}{2j+3},\qquad j\ge1.
\end{equation*}
Using these coefficients and the Legendre multipliers $s_j$,
we obtain the coefficient of $L_j$ in $(S_{n,k}v)'$ as
\begin{equation*}
 s_{j+1}b_j+(2j+1)
 \sum_{\substack{\ell\ge j+2\\\ell-j\text{ even}}}
 \frac{(s_{\ell+1}-s_{\ell-1})b_\ell}{2\ell+1}.
\end{equation*}
For the first term, using $\norm{v'}_I^2=\sum_{\ell\ge0}2|b_\ell|^2/(2\ell+1)$ and $0\le s_{j+1}\le1$ gives
\begin{equation*}
 \left\|\sum_{j\ge0}s_{j+1}b_jL_j\right\|_I^2
 =\sum_{j\ge0}\frac{2s_{j+1}^2|b_j|^2}{2j+1}
 \le\norm{v'}_I^2.
\end{equation*}
For the second term, since $s_{\ell+1}-s_{\ell-1}=-(\omega_{\ell-1}+\omega_\ell)$, we have
\begin{equation*}
 |s_{\ell+1}-s_{\ell-1}|\le\frac{2\ell+1}{Z},\qquad
 s_{\ell+1}-s_{\ell-1}=0\quad\text{if }\ell<n\text{ or }\ell>k+1.
\end{equation*}
Cauchy--Schwarz gives
\begin{equation*}
 \begin{aligned}
 &\left\|\sum_{j\ge0}(2j+1)L_j
 \sum_{\substack{\ell\ge j+2\\\ell-j\text{ even}}}
 \frac{(s_{\ell+1}-s_{\ell-1})b_\ell}{2\ell+1}\right\|_I^2\\
 &=\sum_{j\ge0}2(2j+1)
 \left|\sum_{\substack{\ell\ge j+2\\\ell-j\text{ even}}}
 \frac{(s_{\ell+1}-s_{\ell-1})b_\ell}{2\ell+1}\right|^2\\
 &\le\sum_{j\ge0}2(2j+1)
 \left(\sum_{\substack{\ell\ge j+2\\\ell-j\text{ even}}}
 \frac{|s_{\ell+1}-s_{\ell-1}|^2}{2\ell+1}\right)
 \left(\sum_{\substack{\ell\ge j+2\\\ell-j\text{ even}}}
 \frac{|b_\ell|^2}{2\ell+1}\right)\\
 &\le\left(\sum_{\ell=n}^{k+1}
 \frac{|s_{\ell+1}-s_{\ell-1}|^2}{2\ell+1}
 \sum_{\substack{0\le j\le\ell-2\\\ell-j\text{ even}}}(2j+1)\right)
 \norm{v'}_I^2\\
 &\le\frac1{2Z^2}\sum_{\ell=n}^{k+1}(2\ell+1)\ell(\ell-1)
 \norm{v'}_I^2
 \le\frac{2k(k+1)}Z\norm{v'}_I^2.
 \end{aligned}
\end{equation*}
Here we used $\sum_{\substack{0\le j\le\ell-2\\\ell-j\text{ even}}}(2j+1)
 =\frac{\ell(\ell-1)}2,\qquad
 \sum_{\ell=n}^{k+1}(2\ell+1)
 =\frac{2Z(k-n+2)}{k-n+1}\le4Z,$
and $\ell(\ell-1)\le k(k+1)$ for $n\le\ell\le k+1$.
Taking square roots and applying the triangle inequality proves the
derivative estimate.
\endgroup

Using $L_j(1)=1$, we obtain
\begin{equation*}
 (v-P_mv)(1)
 =\sum_{j=m+1}^{\infty}c_j
 =\frac{b_m}{2m+1}+\frac{b_{m+1}}{2m+3}.
\end{equation*}

Since $\sum_{m=n}^k\omega_m=1$ and $\omega_{n-1}=\omega_{k+1}=0$, we have
\begin{equation*}
 \begin{aligned}
 (v-S_{n,k}v)(1)
 &=\sum_{m=n}^k\omega_m(v-P_mv)(1)\\
 &=\sum_{m=n}^k\omega_m
 \left(\frac{b_m}{2m+1}+\frac{b_{m+1}}{2m+3}\right)\\
 &=\sum_{\ell=n}^k\frac{\omega_\ell b_\ell}{2\ell+1}
 +\sum_{\ell=n+1}^{k+1}\frac{\omega_{\ell-1}b_\ell}{2\ell+1}\\
 &=\sum_{\ell=n}^{k+1}
 \frac{\omega_{\ell-1}+\omega_\ell}{2\ell+1}\,b_\ell.
 \end{aligned}
\end{equation*}
Using the preceding bound and weighted Cauchy--Schwarz, we obtain
\begin{equation*}
 \begin{aligned}
 |(v-S_{n,k}v)(1)|^2
 &\le\left(\sum_{\ell=n}^{k+1}\frac{2\ell+1}{2}{1\over Z^2}\right)
 \left(\sum_{\ell=n}^{k+1}\frac{2|b_\ell|^2}{2\ell+1}\right)\\
 &\le\frac1{2Z^2}\sum_{\ell=n}^{k+1}(2\ell+1)\norm{v'}_I^2\\
 &\le\frac2Z\norm{v'}_I^2.
 \end{aligned}
\end{equation*}

Alternating signs at $-1$ give the same bound.
Polynomial density in $H^1(I)$ proves both bounds for general $v$.
\end{proof}

\subsection{Endpoint correction preserving the moments}

Let $\varphi_j=L_j/\norm{L_j}_I$, $j\ge0$, and
$H=\Pp_r\ominus_{L^2}\Pp_k
=\operatorname{span}\{\varphi_j\}_{j=k+1}^r$.
For $\delta=(\delta_-,\delta_+)^\top\in\R^2$, let $h_\delta$ be the
minimum-$L^2$-norm member of $H$ such that
\begin{equation*}
 \begin{cases}
 h_\delta(-1)=\delta_-,\\
 h_\delta(1)=\delta_+.
 \end{cases}
\end{equation*}
Let $q_\pm(x)=\sum_{j=k+1}^r\varphi_j(\pm1)\varphi_j(x)$,
then $h_\delta(\pm1)=(h_\delta,q_\pm)_I$.
By minimality, $h_\delta=aq_-+bq_+$, and the endpoint conditions give
\begin{equation*}
 \underbrace{\begin{pmatrix}
 (q_-,q_-)_I&(q_-,q_+)_I\\
 (q_+,q_-)_I&(q_+,q_+)_I
 \end{pmatrix}}_{G}
 \begin{pmatrix}a\\b\end{pmatrix}
 =\begin{pmatrix}\delta_-\\\delta_+\end{pmatrix}.
\end{equation*}
A direct computation shows that the smallest eigenvalue of $G$ satisfies
$\gamma\ge(r-k)(r+k+2)/4$. Hence
\begin{equation*}
 \norm{h_\delta}_I^2
 =\int_I(G^{-1}\delta)^{\top}
 (q_-,q_+)^{\top}(q_-,q_+)G^{-1}\delta
 =\delta^\top G^{-1}\delta\le|\delta|^2\gamma^{-1}.
\end{equation*}
We have $L_j'=\sum_{\substack{0\le i<j\\j-i\text{ odd}}}(2i+1)L_i$, then $\varphi_j'=\sum_{\substack{0\le i<j\\j-i\text{ odd}}}
\sqrt{(2i+1)(2j+1)}\,\varphi_i$.
For $0\le i,j\le r$, set
\begin{equation*}
 D_{ij}=\begin{cases}
 \sqrt{(2i+1)(2j+1)},&i<j,\quad j-i\text{ odd},\\
 0,&\text{otherwise}.
 \end{cases}
\end{equation*}
The derivative matrix $D$ gives
\begin{equation*}
 \norm{h_\delta'}_I^2
 \le\left(\sum_{i,j=0}^r|D_{ij}|^2\right)\norm{h_\delta}_I^2
 \le\frac{(r+1)^4}{4}\norm{h_\delta}_I^2
 \le C^2(r+1)^4\gamma^{-1}|\delta|^2.
\end{equation*}
From the preceding derivation, we obtain
\begin{equation}\label{eq:repair}
 \norm{h_\delta}_I\le\gamma^{-1/2}|\delta|,
 \qquad
 \norm{h_\delta'}_I\le C(r+1)^2\gamma^{-1/2}|\delta|.
\end{equation}

Fix a reference triangle $\widehat T$. Let $\Gamma$ be one edge of
$\partial \widehat T$, the open chain formed by two adjacent edges of $\partial \widehat T$,
or the whole boundary $\partial \widehat T$. For either open chain, impose zero
boundary conditions at its two outer endpoints.

For each edge $e$, apply
$S_{n,k}$ in an affine coordinate and denote the resulting polynomial by
$\widetilde g_e$. At a vertex $z$, prescribe
\begin{equation*}
 c_z=\begin{cases}
 0,&z\text{ is an outer endpoint},\\
 \bigl(\widetilde g_{e_1}(z)+\widetilde g_{e_2}(z)\bigr)/2,
   &z=e_1\cap e_2\text{ is a shared vertex}.
 \end{cases}
\end{equation*}
If $e$ has endpoints $z_e^-$ and $z_e^+$, set
\begin{equation*}
 \delta_e=\bigl(c_{z_e^-}-\widetilde g_e(z_e^-),
                 c_{z_e^+}-\widetilde g_e(z_e^+)\bigr)^{\top}.
\end{equation*}
Let $h_e$ be the affine image on $e$ of the correction $h_{\delta_e}$
defined above. We define $\Pi_{n,r}$ edgewise by
\begin{equation*}
 (\Pi_{n,r}g)|_e=\widetilde g_e+h_e.
\end{equation*}
In the operator estimates below, the $H^1$ input on an open chain is
required
to vanish at the two exterior endpoints.
On the closed boundary, we use $H^{1/2}$, and $\Pi_{n,r}$
preserves constants. Every edge integral is preserved in
either setting. All constants depend only on \mbox{reference geometry}.

\begin{proposition}\label{prop:pi}
The map $\Pi_{n,r}$ acts on $L^2(\Gamma)$ and has continuous
edgewise-$\Pp_r$ output that vanishes at the two exterior endpoints when
$\Gamma$ is open
and preserves the \mbox{edge moments}
\begin{equation}\label{eq:pi-mom}
 \int_e(\Pi_{n,r}g-g)\zeta=0,\qquad\zeta\in\Pp_n(e).
\end{equation}
Set $\kappa=(r+1)/(r-n)$.
The operator satisfies
\begin{equation*}
 \norm{\Pi_{n,r}}_{L^2\to L^2}\le C\kappa^{1/2},\enspace
 \norm{\Pi_{n,r}}_{H^1\to H^1}\le C\kappa,\enspace
 \norm{\Pi_{n,r}}_{H^{1/2}_{00}\to H^{1/2}_{00}}\le C\kappa^{3/4}.
\end{equation*}
\end{proposition}
\begin{proof}
\begin{equation*}
 \int_e(\Pi_{n,r}g-g)\zeta
 =\underbrace{\int_e(\widetilde g_e-g)\zeta}_{0}
 +\underbrace{\int_e h_e\zeta}_{0}
 =0,\qquad\zeta\in\Pp_n(e).
\end{equation*}
The first term vanishes by the moment preservation of $S_{n,k}$.
The second vanishes because $h_\delta\in H$ is orthogonal to the affine
pullback $\zeta\in\Pp_n(I)\subset\Pp_k(I)$.

For $m\ge0$ and $p\in\Pp_m(I)$, Cauchy--Schwarz yields
\begin{equation}\label{eq:polynomial-endpoint}
 |p(\pm1)|^2
 \le\left(\sum_{j=0}^m\frac{2j+1}{2}\right)\norm p_I^2
 =\frac{(m+1)^2}{2}\norm p_I^2.
\end{equation}

For an $L^2$ input,
write $\widetilde g|_e=\widetilde g_e$ and $h|_e=h_e$,
with derivatives taken edgewise.
 \eqref{eq:polynomial-endpoint} bounds $\lvert h_e(z)\rvert$ by $C(k+1)$ times
neighboring $L^2$ norms, for each endpoint $z$ of $e$.
By \eqref{eq:repair}, square summation over the edges gives
\begin{equation*}
 \norm h_{L^2}
 \le \frac{C}{\sqrt\gamma}
      \left(\sum_e|\delta_e|^2\right)^{1/2}
 \le \frac{C(k+1)}{\sqrt\gamma}\norm g_{L^2}.
\end{equation*}
For a continuous $H^1$ input the true endpoint
values agree, or are zero at the
two exterior endpoints of an open chain. Thus Lemma~\ref{lem:sn}
bounds $\lvert h_e(z)\rvert$ by $CZ^{-1/2}$ times neighboring derivative norms.
Using \eqref{eq:repair} again, we obtain
\begin{equation*}
 \norm{h'}_{L^2}
 \le \frac{C(r+1)^2}{\sqrt\gamma}
      \left(\sum_e|\delta_e|^2\right)^{1/2}
 \le \frac{C(r+1)^2}{\sqrt{Z\gamma}}\norm{g'}_{L^2}.
\end{equation*}
Lemma~\ref{lem:sn} and \eqref{eq:repair} therefore imply
\begin{align*}
 \norm{\Pi_{n,r}g}_{L^2}
 &\le \norm{\widetilde g}_{L^2}+\norm h_{L^2}
 \le C\left(1+\frac{k+1}{\sqrt\gamma}\right)\norm g_{L^2},\\
 \norm{(\Pi_{n,r}g)'}_{L^2}
 &\le \norm{\widetilde g'}_{L^2}+\norm{h'}_{L^2}
 \le C\left(1+\frac{k+2}{\sqrt Z}
 +\frac{(r+1)^2}{\sqrt{Z\gamma}}\right)\norm{g'}_{L^2}.
\end{align*}
The parameter choice gives
\begin{equation*}
 k-n+1\ge\frac{r-n}{3},\quad r-k\ge\frac{r-n}{2},\quad
 Z\ge\frac{(r-n)(r+1)}{12},\quad
 \gamma\ge\frac{(r-n)(r+1)}8.
\end{equation*}
The bounds on $Z$ and $\gamma$ give the $L^2$ and $H^1$ operator
estimates with factors $C\kappa^{1/2}$ and $C\kappa$, respectively.
Interpolating by \eqref{eq:trace-interpolation} gives the factor
$C\kappa^{3/4}$.
\end{proof}

We abbreviate the trace stability factor by
\begin{equation*}
 \mathfrak T_{n,r}=\left(\frac{r+1}{r-n}\right)^{3/4},
 \qquad n\ge1,\quad r\ge n+2.
\end{equation*}

\section[Moment preservation on element bubbles]{%
  \texorpdfstring{A moment-preserving operator on $H_0^1(\widehat T)$}%
  {Moment preservation on element bubbles}}\label{sec:bubble}
In this and the next section we work on the reference triangle
\begin{equation*}
 \begin{gathered}
 \widehat T=\{(x,y):x>0,\ y>0,\ x+y<1\},\qquad b=xy(1-x-y),\\
 a_{\widehat T}(v,w)=(\nabla v,\nabla w)_{\widehat T}.
 \end{gathered}
\end{equation*}
We define the weighted and unweighted exact-degree spaces by
\begin{equation*}
 Q_j^b=\Pp_j\ominus_{L^2_b}\Pp_{j-1},\quad
 Q_j^0=\Pp_j\ominus_{L^2}\Pp_{j-1},\quad
 W_j=bQ_j^b,\quad (z,t)_{L^2_b}=\int_{\widehat T} bzt.
\end{equation*}
The union $\bigcup_{m\ge0}b\Pp_m(\widehat T)$ is dense in $H^1_0(\widehat T)$. Indeed, for
$w\in C_c^\infty(\widehat T)$ we approximate $w/b$ in
$C^1(\overline{\widehat T})$ by
polynomials and multiply by $b$. Density of $C_c^\infty(\widehat T)$ gives
the assertion.
For $w_j=bz_j\in W_j$ and $w_k\in W_k$, integration by parts
gives $a_{\widehat T}(w_j,w_k)=-(z_j,\Delta w_k)_{L^2_b}$ with
$\Delta w_k\in\Pp_{k+1}$. The defining orthogonality of $Q_j^b$ and
symmetry give the first identity below. The second follows from
$(F,w_j)_{\widehat T}=(F,z_j)_{L^2_b}$.
\begin{equation}\label{eq:band}
 a_{\widehat T}(W_j,W_k)=0\ (|j-k|>1),\qquad
 (F,W_j)_{\widehat T}=0\ (F\in\Pp_\ell,\ j>\ell).
\end{equation}

\subsection{An angle bound for neighboring bubble subspaces}
\begin{lemma}\label{lem:split}
For $w_j\in W_j$ and $w_{j+1}\in W_{j+1}$, we have
\begin{equation}\label{eq:angle}
 |a_{\widehat T}(w_j,w_{j+1})|\le\gamma_j\norm{\nabla w_j}_{\widehat T}
 \norm{\nabla w_{j+1}}_{\widehat T},\quad
 \gamma_j^2=\frac{(j+2)(j+4)}{(2j+5)(2j+7)}<\frac14.
\end{equation}
\end{lemma}
\begin{proof}
Let $w=bz\in W_j$. We first establish the orthogonal decomposition
$\nabla w=g^-+g^+$ with $g^-\in[Q_{j+1}^0]^2$ and
$g^+\in[Q_{j+2}^0]^2$, and the identities
\begin{equation*}
 \norm{g^-}_{\widehat T}^2=\frac{j+3}{2j+5}\norm{\nabla w}_{\widehat T}^2,
 \qquad \norm{g^+}_{\widehat T}^2=\frac{j+2}{2j+5}\norm{\nabla w}_{\widehat T}^2.
\end{equation*}

For $\boldsymbol q\in[\Pp_j]^2$, integration by parts gives
$(\nabla w,\boldsymbol q)_{\widehat T}=-(bz,\diver\boldsymbol q)_{\widehat T}=0$.
Thus we have
\begin{equation*}
 \nabla w\in[\Pp_{j+2}]^2\ominus_{L^2}[\Pp_j]^2
 =[Q_{j+1}^0]^2\oplus[Q_{j+2}^0]^2.
\end{equation*}

To compute the $L^2(\widehat T)$ norms of $g^-$ and $g^+$, let
$e_y=\{(x,0):0<x<1\}$. For $m\ge0$ and $p\in Q_m^0$, take $V=(x,y-1)$.
Then $\diver V=2$ and $V\cdot\nabla p-mp\in\Pp_{m-1}$, while $V\cdot n=1$ on $e_y$
and $V\cdot n=0$ on the other two edges.
Orthogonality gives
\begin{equation*}
 (p,V\cdot\nabla p-mp)_{\widehat T}=0,\qquad
 (p,V\cdot\nabla p)_{\widehat T}=m\norm p_{\widehat T}^2.
\end{equation*}
The divergence theorem then gives
\begin{equation}\label{eq:exact-degree-trace}
 \norm p_{e_y}^2
 =2\norm p_{\widehat T}^2+2(p,V\cdot\nabla p)_{\widehat T}
 =2(m+1)\norm p_{\widehat T}^2.
\end{equation}
Write $g^\pm=(g_1^\pm,g_2^\pm)$. Since $w_x=0$ on $e_y$,
we have $g_1^-=-g_1^+$ there. Applying \eqref{eq:exact-degree-trace} with
$m=j+1$ and $m=j+2$ yields
\begin{equation*}
 2(j+2)\norm{g_1^-}_{\widehat T}^2
 =\norm{g_1^-}_{e_y}^2
 =\norm{g_1^+}_{e_y}^2
 =2(j+3)\norm{g_1^+}_{\widehat T}^2.
\end{equation*}
Interchanging $x$ and $y$ and using $w_y=0$ on $x=0$ gives
\begin{equation*}
 (j+2)\norm{g_2^-}_{\widehat T}^2=(j+3)\norm{g_2^+}_{\widehat T}^2.
\end{equation*}
Writing $X_\pm=\norm{g^\pm}_{\widehat T}^2$, orthogonality and addition of the
component identities give
\begin{equation*}
 X_-+X_+=\norm{\nabla w}_{\widehat T}^2,\qquad
 (j+2)X_-=(j+3)X_+.
\end{equation*}

Solving these equations gives the two squared-norm identities.
Let $g^\pm(w_k)$ be the two components of $\nabla w_k$ for
$k=j,j+1$. Orthogonality of the different degree spaces and Cauchy--Schwarz
give
\begin{equation*}
 |a_{\widehat T}(w_j,w_{j+1})|
 = |(g^+(w_j),g^-(w_{j+1}))_{\widehat T}|
 \le \norm{g^+(w_j)}_{\widehat T}\norm{g^-(w_{j+1})}_{\widehat T}.
\end{equation*}
Substituting the squared-norm identities gives \eqref{eq:angle}.
\end{proof}

\subsection{A stable cosine cutoff in polynomial degree}
For a bubble polynomial $v=\sum_{j=0}^J v_j$, with $v_j\in W_j$, set
$a_j=\norm{\nabla v_j}_{\widehat T}$. Set $v_j=0$ and $a_j=0$ for $j>J$.
Combining \eqref{eq:band} and
\eqref{eq:angle} gives
\begin{equation*}
  \norm{\nabla v}_{\widehat T}^2
 \ge\sum_{j=0}^J a_j^2-\sum_{j=0}^{J-1}a_ja_{j+1}
 =\frac12\left(a_0^2+a_J^2+
 \sum_{j=0}^{J-1}(a_{j+1}-a_j)^2\right).
\end{equation*}
It follows from the identity $a_j=a_0+\sum_{i=0}^{j-1}(a_{i+1}-a_i)$
and Cauchy--Schwarz that
\begin{equation}\label{eq:level-bound}
 a_j\le\sqrt{2(j+1)}\norm{\nabla v}_{\widehat T}.
\end{equation}

For integers $\ell\ge0$ and $r\ge\ell+3$, let $L=r-\ell-2$ and define
\begin{equation*}
 \theta_j=\begin{cases}
 0,&j\le\ell,\\
 \pi(j-\ell)/(2L),&\ell<j<\ell+L,\\
 \pi/2,&j\ge\ell+L,
 \end{cases}
 \qquad c_j=\cos\theta_j,\quad s_j=\sin\theta_j.
\end{equation*}
We set $B_{\ell,r}v=\sum_jc_jv_j$. Since $c_j=0$ for $j\ge r-2$,
the output belongs to $b\Pp_{r-3}\subset\Pp_r$.

\begin{proposition}\label{prop:bubble}
The map $B_{\ell,r}$ extends to a bounded linear operator
$H^1_0(\widehat T)\to H^1_0(\widehat T)\cap\Pp_r(\widehat T)$ satisfying
\begin{align}
 (F,B_{\ell,r}v)_{\widehat T}&=(F,v)_{\widehat T}\qquad(F\in\Pp_\ell(\widehat T)),\label{eq:b-mom}\\
 \norm{\nabla B_{\ell,r}v}_{\widehat T}&\le\mathfrak B_{\ell,r}\norm{\nabla v}_{\widehat T},
 \qquad
 \mathfrak B_{\ell,r}^2=1+\frac{\pi^2(\ell+r+1)}{8(r-\ell-2)}.\notag
\end{align}
\end{proposition}
\begin{proof}
It follows from \eqref{eq:band} and $c_j=1$ for $j\le\ell$ that
\begin{equation*}
 (F,B_{\ell,r}v-v)_{\widehat T}
 =\sum_{j=0}^{J}(c_j-1)(F,v_j)_{\widehat T}
 =0,\qquad F\in\Pp_\ell(\widehat T).
\end{equation*}
It remains to prove the energy bound. Let $D_sv=\sum_js_jv_j$.
The tridiagonal identity \eqref{eq:band} and $c_j^2+s_j^2=1$ give
\begin{equation}\label{eq:cutoff-energy}
 \begin{aligned}
 &\norm{\nabla B_{\ell,r}v}_{\widehat T}^2+\norm{\nabla D_sv}_{\widehat T}^2
 =\norm{\nabla v}_{\widehat T}^2\\
 &\qquad-\sum_{j\ge0}\bigl((c_{j+1}-c_j)^2+(s_{j+1}-s_j)^2\bigr)
 a_{\widehat T}(v_j,v_{j+1}).
 \end{aligned}
\end{equation}
For the second term, only $\ell\le j<\ell+L$ contribute. So we have
\begin{equation*}
 (c_{j+1}-c_j)^2+(s_{j+1}-s_j)^2
 =4\sin^2\!\left(\frac{\pi}{4L}\right)
 \le\frac{\pi^2}{4L^2}.
\end{equation*}
Combining \eqref{eq:angle} and \eqref{eq:level-bound} gives
\begin{equation*}
 |a_{\widehat T}(v_j,v_{j+1})|\le(j+2)\norm{\nabla v}_{\widehat T}^2,
 \qquad \sum_{j=\ell}^{\ell+L-1}(j+2)=\frac{L(\ell+r+1)}2.
\end{equation*}
Substituting these estimates into \eqref{eq:cutoff-energy} and
dropping $\norm{\nabla D_sv}_{\widehat T}^2$
proves the~bound.

By density, take bubble polynomials $v_m\to v$ in $H_0^1(\widehat T)$.
The energy bound applied to $v_m-v_n$ makes $B_{\ell,r}v_m$ Cauchy in the closed space
$b\Pp_{r-3}(\widehat T)$. The same bound makes the limit independent of the
approximating sequence. We define $B_{\ell,r}v$ by this limit and pass the
energy bound to the limit. Both sides of \eqref{eq:b-mom} converge by
Poincar\'e's inequality.
\end{proof}

\subsection{Comparison with the square-domain proof}\label{sec:comparison}
On the square, the tensor-product Babu\v ska--Shen basis is
grouped by total degree \cite[Section~2.1, (2.11)]{CNSV19}.
Distinct degree spaces are energy-orthogonal unless their indices differ
by two \cite[(2.13)]{CNSV19}.
For each total-degree index $j\ge6$, the energy-orthogonal
projection from the degree-$(j-2)$ space onto the degree-$j$ space
has norm at most $\tfrac12(1-c_{\mathrm{sq}}j^{-2})$, with a constant
$c_{\mathrm{sq}}>0$ independent of $j$
\cite[Proposition~2.2]{CNSV19}.
For $q>p+4$, the energy norm of each of the two
highest-degree components is at most $6/(q-p)$ times that of the
square Galerkin solution \cite[Proposition~3.1]{CNSV19}.
The test-function tail above degree $q$ has energy norm at most $Cq$
times the energy norm of the test function \cite[Proposition~4.2]{CNSV19}.
Combining these estimates gives the factor $q/(q-p)$ in the saturation
bound \cite[Section~5]{CNSV19}.

The analysis relies on the Cartesian tensor-product structure of the
square \cite[p.~990]{CNSV19} and therefore does not apply directly to the
triangle. Even after pullback to the square by the Duffy map, the
Dirichlet energy contains weights and cross terms, and the polynomial
test spaces differ from those used in the square proof.
We establish the required energy-angle estimate and norm comparison
directly on the triangle through the gradient decomposition of bubble
functions and the cosine cutoff.

On the reference triangle, Lemma~\ref{lem:split} splits the
gradient of each bubble component into two orthogonal polynomial
components of consecutive degrees and computes their exact energy
fractions. These fractions give the explicit angle bound
\eqref{eq:angle} for consecutive bubble levels. We then apply a cosine
cutoff in the degree index to arbitrary $H^1_0(\widehat T)$ test functions.
In \eqref{eq:cutoff-energy}, the squared increments are $O(L^{-2})$
over $L=r-\ell-2$ levels. Combining this identity with
\eqref{eq:level-bound} gives the factor $\mathfrak B_{\ell,r}$.
Moment preservation and Lemma~\ref{lem:fortin} consequently yield
\begin{equation*}
 \dual F{H^1_0(\widehat T)}\le\mathfrak B_{\ell,r}
 \dual F{H^1_0(\widehat T)\cap\Pp_r(\widehat T)},\qquad
 F\in\Pp_\ell(\widehat T),\quad r\ge\ell+3.
\end{equation*}
The cosine cutoff yields the volume dual norm comparison
with a stability factor bounded by $C(r/(r-\ell-2))^{1/2}$.
The boundary estimates use the moment-preserving trace operator of
Proposition~\ref{prop:pi}, and the mixed estimate also uses the polynomial
adjoint identity of Lemma~\ref{lem:adjoint}. Polynomial flux liftings
carry these local estimates to vertex patches, where they bound the
continuous residual dual norm by its discrete counterpart uniformly in
the polynomial degree for each fixed enrichment factor $\lambda>1$.

\section{Trace liftings and reference estimates}\label{sec:mixed}
\subsection{Degree-preserving trace liftings}
We use variants of the averaging liftings in
\cite[Proposition~B.1, (B.2) and (B.8)]{M03} and prove the estimates needed
for the mixed problem. For smooth boundary data
$g:(0,1)\to\R$ on the edge $y=0$, we define, for $(x,y)\in \widehat T$,
\begin{equation}\label{eq:e2}
\begin{aligned}
 (\mathcal E_0g)(x,y)&=\frac1y\int_x^{x+y}g(s)\,ds,\\
 (\mathcal E_1g)(x,y)&=\frac{x(1-x-y)}y
 \int_x^{x+y}\frac{g(s)}{s(1-s)}\,ds.
\end{aligned}
\end{equation}
Both $\mathcal E_0g$ and $\mathcal E_1g$ have trace $g$ on $y=0$, and $\mathcal E_1g$
vanishes on the other two edges.

Subtracting $\mathcal E_1g$ from a
mixed test function with free-edge trace $g$ leaves a remainder in
$H^1_0(\widehat T)$. The Schwarz decomposition
in~\cite[Sections~3.3--3.4]{SMPZ08} likewise subtracts
lifted edge traces after removing the vertex contributions, leaving an
element bubble. Both constructions require a degree-preserving lifting
with an energy bound independent of the degree, and the edge lifting used there
agrees with $\mathcal E_1$ up to an affine change of coordinates. In
Theorem~\ref{thm:mixed}, we use its adjoint to express a weighted volume
pairing as a polynomial edge moment.

\begin{lemma}\label{lem:lifts}
The operator $\mathcal E_0$ is bounded from $H^{1/2}(0,1)$ to $H^1(\widehat T)$.
The operator $\mathcal E_1$ is bounded from $H^{1/2}_{00}(0,1)$ to $H^1(\widehat T)$.
For polynomial data $g\in\Pp_r(0,1)$, we have $\mathcal E_0g\in\Pp_r(\widehat T)$.
If $g(0)=g(1)=0$, then $\mathcal E_1g\in\Pp_r(\widehat T)$.
\end{lemma}
\begin{proof}
\emph{Estimate for $\mathcal E_0$.}
We use the coordinates $(a,c)=(x,x+y)$ and write $\mathcal E_i g$ for
$(a,c)\mapsto(\mathcal E_i g)(a,c-a)$.
We set $s_t=(1-t)a+tc$.
By Jensen's inequality and Fubini's theorem,
\begin{equation*}
 \norm{\mathcal E_0g}_{L^2(\widehat T)}^2
 \le \int_0^1\int_a^1\frac{1}{c-a}\int_a^c |g(s)|^2\,ds\,dc\,da
 =\int_0^1 |g(s)|^2
     \left(\int_0^s\int_s^1\frac{dc\,da}{c-a}\right)ds .
\end{equation*}
The weight satisfies
\begin{equation*}
 \int_0^s\int_s^1\frac{dc\,da}{c-a}
 =-s\log s-(1-s)\log(1-s)\le\log 2,
\end{equation*}
so we obtain
$\norm{\mathcal E_0g}_{L^2(\widehat T)}^2\le(\log 2)\norm{g}_{L^2(0,1)}^2$.
Also
\begin{equation*}
 \partial_a\mathcal E_0g=\frac1{(c-a)^2}\int_a^c(g(s)-g(a))\,ds.
\end{equation*}
Cauchy--Schwarz followed by integration in $c$ yields
\begin{equation*}
 \int_{0<a<c<1}|\partial_a\mathcal E_0g|^2\,da\,dc
 \le\frac12\int_{0<a<s<1}\frac{|g(s)-g(a)|^2}{(s-a)^2}\,da\,ds.
\end{equation*}
Reflecting $(a,c,s)$ to $(1-c,1-a,1-s)$ gives the bound
for $\partial_c \mathcal E_0g$. Thus we have
$\norm{\mathcal E_0g}_{H^1(\widehat T)}\le C\norm{g}_{H^{1/2}(0,1)}$.

\emph{Weighted estimates for $\mathcal E_1$.}
For any fixed $0<t<1$, changing variables from $(a,c)$ to $(s_t,c-a)$ gives
\begin{equation}\label{eq:weighted-change}
 \norm{\frac{g(s_t)}{s_t}}_{L^2(a<c)}^2
 \le\frac1t\int_0^1\frac{|g(s)|^2}s\,ds,
 \quad
 \norm{\frac{g(s_t)}{1-s_t}}_{L^2(a<c)}^2
 \le\frac1{1-t}\int_0^1\frac{|g(s)|^2}{1-s}\,ds.
\end{equation}
Indeed the allowed interval of $c-a$ has length at most $s_t/t$, and also
at most $(1-s_t)/(1-t)$.
For $F_1=\mathcal E_1g-\mathcal E_0g$, direct differentiation gives
\begin{align*}
 (F_1)_a={}&(1-c)\int_0^1(1-t)\frac{g(s_t)}{s_t}\,dt
 -c\int_0^1(1-t)\frac{g(s_t)}{1-s_t}\,dt+\frac{g(a)}{1-a},\\
 (F_1)_c={}&(1-a)\int_0^1t\frac{g(s_t)}{s_t}\,dt
 -a\int_0^1t\frac{g(s_t)}{1-s_t}\,dt-\frac{g(c)}c.
\end{align*}
We estimate the integral terms using Minkowski's inequality and
\eqref{eq:weighted-change}. We integrate the endpoint terms directly. Hence
\begin{equation*}
 \norm{\nabla F_1}_{\widehat T}^2
 \le C\int_0^1\frac{|g(s)|^2}{s(1-s)}\,ds.
\end{equation*}
Moreover,
$0\le a(1-c)/(s(1-s))\le1$ for $a<s<c$.
Consequently,
\begin{equation*}
 |\mathcal E_1g|\le \mathcal E_0|g|.
\end{equation*}
The $L^2$ bound for $\mathcal E_1g$ follows from that for $\mathcal E_0|g|$.
Together with $\mathcal E_1g=\mathcal E_0g+F_1$, the gradient bounds for $\mathcal E_0g$ and $F_1$
give the $H^1$ bound for $\mathcal E_1g$.

\emph{Degree preservation.}
For $g\in\Pp_r(0,1)$, if $g(0)=g(1)=0$, then
$g(s)/(s(1-s))$ is a polynomial of degree at most $r-2$.
An interval average preserves polynomial degree in $(a,c)$.
The factor $a(1-c)$ in $\mathcal E_1g$ gives total degree at most $r$.
The trace statements hold first for smooth admissible data and
then by density and boundedness in the indicated trace spaces.
\end{proof}

\begin{corollary}
\label{cor:scalar-extension}
Let $\Gamma$ be one, two, or three edges of $\widehat T$. If compatible edgewise
$\Pp_r$ data $g$ on $\Gamma$ are the trace of an $H^1(\widehat T)$ function, they
have a polynomial extension $L_\Gamma g\in\Pp_r(\widehat T)$ satisfying
\begin{equation*}
 \norm{L_\Gamma g}_{H^1(\widehat T)}\le C\inf_{\tr_\Gamma v=g}
 \norm v_{H^1(\widehat T)}.
\end{equation*}
There is also the seminorm version, with the infimum of
$\norm{\nabla v}_{\widehat T}$, and a suitable choice of the extension modulo constants.
\end{corollary}
\begin{proof}
Let $L_\Gamma$ be the linear lifting from
\cite[Theorem~B.4 and Remark~B.5]{M03}, with $q=2$ and transported to $\widehat T$
by a fixed affine map. It preserves constants and polynomial degree.
For any $v\in H^1(\widehat T)$ with $\tr_\Gamma v=g$, the lifting bound and the
trace theorem give
\begin{equation*}
 \norm{L_\Gamma g}_{H^1(\widehat T)}
 \le C\norm g_{H^{1/2}(\Gamma)}
 \le C\norm v_{H^1(\widehat T)}.
\end{equation*}
With $\bar v=|\widehat T|^{-1}\int_{\widehat T} v$, constant preservation and Poincar\'e's
inequality give
\begin{equation*}
 \norm{\nabla L_\Gamma g}_{\widehat T}
 =\norm{\nabla L_\Gamma(g-\bar v)}_{\widehat T}
 \le C\norm{v-\bar v}_{H^1(\widehat T)}
 \le C\norm{\nabla v}_{\widehat T}.
\end{equation*}
Taking the infimum over $v$ yields both bounds.
\end{proof}

\subsection{Dual norm estimates for volume loads with mixed boundary conditions}
Let $X=H^1_{\{x=0\}\cup\{x+y=1\}}(\widehat T)$ and $F=y\zeta$, with
$\zeta\in\Pp_{\ell-1}(\widehat T)$ and $\ell\ge1$. Here $y$ is the barycentric coordinate
equal to one at the intersection of the two essential edges and zero on the
free edge. For $g\in H^{1/2}_{00}(0,1)$, Fubini gives
\begin{equation}\label{eq:adjoint}
      \begin{aligned}
             (F,\mathcal E_1g)_{\widehat T}&=\int_0^1\eta_F(s)g(s)\,ds,\\
 \eta_F(s)&=\frac1{s(1-s)}\int_0^s\int_s^1
 x(1-t)\zeta(x,t-x)\,dt\,dx.
      \end{aligned}
\end{equation}
Here the factor $y$ in $F$ cancels the interval length in
\eqref{eq:e2}. This cancellation yields a polynomial boundary density.

\begin{lemma}\label{lem:adjoint}
The function $\eta_F$ in \eqref{eq:adjoint} is a polynomial of degree at
most $\ell+1$, divisible by $s(1-s)$.
\end{lemma}
\begin{proof}
$\int_0^s\int_s^1 x(1-t)\zeta(x,t-x)\,dt\,dx$
has degree at most $\ell+3$. We expand
$\zeta(x,t-x)$ into monomials $x^it^j$, with $i+j\le \ell-1$.
The $x$ integral of each term has a factor $s^{i+2}$, whereas the integral of
$(1-t)t^j$ from $s$ to $1$ has a double zero at $s=1$.
The numerator is therefore divisible by $s^2(1-s)^2$.
Division proves both assertions. Formula~\eqref{eq:adjoint} holds first
for $g\in C_c^\infty(0,1)$. The lifting estimate and the polynomial
expression for $\eta_F$ extend it to $H^{1/2}_{00}(0,1)$.
\end{proof}

\begin{theorem}\label{thm:mixed}
Let $X=H^1_{\{x=0\}\cup\{x+y=1\}}(\widehat T)$. For $F=y\zeta$,
$\zeta\in\Pp_{\ell-1}(\widehat T)$, $\ell\ge1$, and $r\ge\ell+3$,
\begin{equation*}
 \dual F X\le C\bigl(\mathfrak B_{\ell,r}+\mathfrak T_{\ell+1,r}\bigr)
 \dual F{(X\cap\Pp_r(\widehat T))}.
\end{equation*}
\end{theorem}
\begin{proof}
For $v\in X$, set $g=\tr_{y=0}v$, $w=v-\mathcal E_1g\in H^1_0(\widehat T)$, and
\begin{equation*}
 J_{\ell,r}v=B_{\ell,r}w+\mathcal E_1\Pi_{\ell+1,r}g.
\end{equation*}
The full trace of $v$ is the zero extension of $g$. The trace
theorem, Poincar\'e's inequality on $X$, and Lemma~\ref{lem:lifts} give
$\norm{\nabla w}_{\widehat T}+\norm g_{H^{1/2}_{00}(0,1)}
\le C\norm{\nabla v}_{\widehat T}$. Combining the bounds for
$B_{\ell,r}$ and $\Pi_{\ell+1,r}$ in Propositions~\ref{prop:bubble}
and~\ref{prop:pi} with the lifting estimate gives
\begin{equation*}
 \norm{\nabla J_{\ell,r}v}_{\widehat T}
 \le C\bigl(\mathfrak B_{\ell,r}+\mathfrak T_{\ell+1,r}\bigr)
 \norm{\nabla v}_{\widehat T}.
\end{equation*}
$B_{\ell,r}w$ and $\mathcal E_1\Pi_{\ell+1,r}g$
have degree at most $r$ and satisfy the essential conditions.
The bubble term preserves the $F$ moment by \eqref{eq:b-mom}.
For the lifted trace, \eqref{eq:adjoint}, Lemma~\ref{lem:adjoint}, and
\eqref{eq:pi-mom} give
\begin{equation*}
 (F,\mathcal E_1\Pi_{\ell+1,r}g)_{\widehat T}=(\eta_F,\Pi_{\ell+1,r}g)_{(0,1)}
 =(\eta_F,g)_{(0,1)}=(F,\mathcal E_1g)_{\widehat T}.
\end{equation*}
Hence $(F,J_{\ell,r}v)_{\widehat T}=(F,v)_{\widehat T}$, and Lemma~\ref{lem:fortin} completes the proof.
\end{proof}

\subsection{Boundary reference estimates}\label{sec:edge}
\begin{theorem}\label{thm:edge}
Let $e_D$ and $e_L$ be distinct edges of $\widehat T$, $X=H^1_{e_D}(\widehat T)$, and
$F(v)=\int_{e_L}\zeta v$ with $\zeta\in\Pp_n(e_L)$, $n\ge1$. For $r\ge n+2$,
\begin{equation*}
 \dual F X\le C\mathfrak T_{n,r}\dual F{(X\cap\Pp_r(\widehat T))}.
\end{equation*}
\end{theorem}
\begin{proof}
Let $\Gamma=\partial \widehat T\setminus e_D$ be the complementary two-edge
chain. For $v\in X$, set $g=\tr_\Gamma v$ and define the full-boundary
polynomial data
\begin{equation*}
 \widehat g|_\Gamma=\Pi_{n,r}g,\qquad \widehat g|_{e_D}=0,
 \qquad J_{n,r}v=L_{\partial \widehat T}\widehat g.
\end{equation*}
Then $J_{n,r}v\in X\cap\Pp_r(\widehat T)$. The seminorm version of
Corollary~\ref{cor:scalar-extension}, Proposition~\ref{prop:pi}, and the
trace theorem give
\begin{equation*}
 \norm{\nabla J_{n,r}v}_{\widehat T}
 \le C\norm{\widehat g}_{H^{1/2}(\partial \widehat T)}
 \le C\mathfrak T_{n,r}\norm{\nabla v}_{\widehat T}.
\end{equation*}
Moreover, \eqref{eq:pi-mom} gives
\begin{equation*}
 F(J_{n,r}v)
 =\int_{e_L}\zeta\,\Pi_{n,r}g
 =\int_{e_L}\zeta g=F(v).
\end{equation*}
Lemma~\ref{lem:fortin} proves the estimate.
\end{proof}

\begin{theorem}\label{thm:neumann}
Let $X=H^1(\widehat T)/\R$ and $X_r=\Pp_r(\widehat T)/\R$,
and let $\zeta$ be edgewise in $\Pp_n$, $n\ge1$, with
$\int_{\partial \widehat T}\zeta=0$. For $F(v)=\int_{\partial \widehat T}\zeta v$ and $r\ge n+2$,
\begin{equation*}
 \dual F X\le C\mathfrak T_{n,r}
 \norm F_{X_r'}.
\end{equation*}
\end{theorem}
\begin{proof}
Choose the mean-zero representative of $v\in X$, set
$g=\tr_{\partial \widehat T}v$, and define
\begin{equation*}
 w=L_{\partial \widehat T}\Pi_{n,r}g,\qquad
 J_{n,r}v=w-|\widehat T|^{-1}
 \int_{\widehat T}w.
\end{equation*}
Then $J_{n,r}v$ is the mean-zero polynomial representative of an element of $X_r$. By
Corollary~\ref{cor:scalar-extension}, Proposition~\ref{prop:pi}, the trace
theorem, and Poincar\'e's~\mbox{inequality},
\begin{equation*}
 \norm{\nabla J_{n,r}v}_{\widehat T}
 =\norm{\nabla w}_{\widehat T}
 \le C\mathfrak T_{n,r}\norm{\nabla v}_{\widehat T}.
\end{equation*}
Since $F(1)=0$, moment preservation on each edge gives
\begin{equation*}
 F(J_{n,r}v)=F(w)
 =\int_{\partial \widehat T}\zeta\,\Pi_{n,r}g
 =\int_{\partial \widehat T}\zeta g=F(v).
\end{equation*}
Lemma~\ref{lem:fortin} proves the estimate.
\end{proof}

Theorems~\ref{thm:mixed}, \ref{thm:edge}, and \ref{thm:neumann} establish
all three dual norm estimates on the reference triangle.
For $r\ge\ell+3$, define
\begin{equation}\label{eq:patch-ratio}
 \rho=\frac{r+1}{r-\ell-2}.
\end{equation}
The stability factors for $B_{\ell,r}$ and $\Pi_{n,r}$ satisfy
$\mathfrak B_{\ell,r}\le C\rho^{1/2}$ and
$\mathfrak T_{\ell+1,r},\mathfrak T_{\ell,r}\le\rho^{3/4}$.
For fixed $\lambda>1$, set
$r_\ell=\lceil\lambda\ell\rceil+2$. Then
\begin{equation*}
 \sup_{\ell\ge1}\max\left\{
 \mathfrak B_{\ell,r_\ell}+\mathfrak T_{\ell+1,r_\ell},
 \mathfrak T_{\ell,r_\ell}\right\}<\infty.
\end{equation*}
Indeed the common local test degree $r_\ell$
satisfies $r_\ell\ge\ell+3$ and
$\rho\le1+4/(\lambda-1)$ when $r=r_\ell$, which bounds both
stability factors.
An affine pullback gives an equivalent constant-matrix energy
norm. Shape regularity makes the equivalence uniform, so the estimates
transfer to physical~triangles.

\section{Dual norm estimates on vertex patches}\label{sec:patch}
\subsection{Polynomial flux liftings}
The following Raviart--Thomas lifting is stated in
\cite[Lemma~6.1]{CNSV17} and follows from the regularized Poincar\'e operators
in~\cite[Proposition~4.1(i) and Section~4.2]{CM}.
Here $H^{-1}(T)$ is the dual of $H^1_0(T)$ equipped with the gradient
norm.

\begin{theorem}\label{thm:cm}
For every $F\in\Pp_\ell(T)$ there is
$\sigma\in\RT_\ell(T)=[\Pp_\ell(T)]^2+(x,y)^\top\Pp_\ell(T)$ such that
\begin{equation}\label{eq:cm}
 \diver\sigma=F,\qquad
 \norm\sigma_T\le C\dual F{H^1_0(T)},
\end{equation}
with $C$ independent of $\ell$. On each edge, $\sigma\cdot n$ has degree
at most $\ell$.
\end{theorem}

Through the Piola transform, the estimate \eqref{eq:cm} also holds on physical
triangles, with a constant depending on shape regularity. To impose
partial normal data, we use a divergence-free lifting.

\begin{lemma}\label{lem:solenoid}
Let $\Gamma$ consist of one or two connected edges of $T$.
If $\tau\in L^2(T)^2$ is divergence-free with an
edgewise-$\Pp_m$ normal trace on $\Gamma$, then some
$\sigma\in[\Pp_m(T)]^2$~\mbox{satisfies}
\begin{equation*}
 \diver\sigma=0,\qquad \sigma\cdot n=\tau\cdot n\text{ on }\Gamma,
 \qquad \norm\sigma_T\le C\norm\tau_T.
\end{equation*}
\end{lemma}
\begin{proof}
We apply \cite[Section~3]{BPS09} with
$g_n=(\tau\cdot n)|_\Gamma$ and use $\tau$ as the comparison field.
This gives $\sigma\in\RT_m(T)$ with
\[
 \diver\sigma=0,\qquad \sigma\cdot n=g_n\text{ on }\Gamma,\qquad
 \norm{\sigma}_T\le C\norm{\tau}_T
\]
with $C$ independent of $m$. Write $\sigma=a+(x,y)^\top q_m$ with
$a\in[\Pp_m(T)]^2$ and $q_m$ homogeneous of degree $m$. The identity
\[
 0=\diver\sigma=\diver a+(m+2)q_m
\]
gives $q_m=0$ by comparing terms of degree $m$. Thus
$\sigma\in[\Pp_m(T)]^2$.
\end{proof}

\subsection{Patch spaces and polynomial test functions}
Index the mesh vertices by $\{a_i\}_{i\in\mathcal I}$.
Let $\Omega_i$ be the union of triangles incident on a mesh vertex $a_i$ and
$\phi_i$ its piecewise affine hat function. At an interior vertex set
$X_i=H^1(\Omega_i)/\R$. At a boundary vertex let $X_i$ consist of
$H^1(\Omega_i)$ functions vanishing on the two boundary radial edges.
We equip $X_i$ with the norm $\norm{\nabla v}_{\Omega_i}$ and set
$h_i=h_{a_i}=\operatorname{diam}\Omega_i$.
Let $\Gamma_{\mathrm{out}}$ be the union of the edges opposite $a_i$ in the triangles of $\Omega_i$.
A radial edge contains $a_i$,
and the outer edge $o_T$ of $T\subset{\Omega_i}$ is opposite $a_i$.
Thus $\Gamma_{\mathrm{out}}=\bigcup_{T\subset\Omega_i}o_T$.
For an integer $\ell\ge1$, consider
\begin{equation}\label{eq:patch-res}
 R(v)=\sum_{T\subset{\Omega_i}}\int_{T}\phi_iG_Tv
       +\sum_{e\text{ internal radial}}\int_e J_{i,e}v,
 \quad G_T\in\Pp_{\ell-1}(T),\quad J_{i,e}\in\Pp_\ell(e).
\end{equation}
For an interior patch assume $R(1)=0$. Let
\begin{equation*}
 X_{i,r}=\{v\in X_i:v|_{T}\in\Pp_r(T)
 \ \forall T\in\T\text{ with }T\subset\Omega_i\}.
\end{equation*}
\begin{lemma}\label{lem:patch-extension}
Let $T_0$ be any triangle in a vertex patch $\Omega_i$, and let $v\in\Pp_r(T_0)$ with $r\ge1$.
On an interior patch, $v$ has a continuous
piecewise-$\Pp_r$ extension $\bar v$ to $\Omega_i$ such~that
\begin{equation*}
 \norm{\nabla\bar v}_{\Omega_i}\le C\norm{\nabla v}_{T_0}.
\end{equation*}
The constant depends only on shape regularity.
On either an interior or a boundary patch, if
a triangle $T_1$ in $\Omega_i$ shares a radial edge with $T_0$
and $v$ vanishes on the other radial edge of $T_0$, then
there is a continuous piecewise-$\Pp_r$ extension $\bar v$
satisfying the same estimate and vanishing
on the other radial edge of $T_1$ and on all triangles
outside $T_0\cup T_1$.
\end{lemma}
\begin{proof}
For the second assertion, let adjacent triangles $T_0$ and $T_1$ share the
radial edge $e$, and let $f$ and $f'$ be their other radial edges. The
affine vertex map $\Phi:T_1\to T_0$ that fixes $e$ pointwise maps $f'$ onto $f$.
If $v|_f=0$, define $\tilde{v}$ on $T_0\cup T_1$ by
\begin{equation*}
 \tilde{v}|_{T_0}=v,\qquad \tilde{v}|_{T_1}=v\circ\Phi.
\end{equation*}
The two pieces agree on $e$, $\tilde{v}|_{f'}=0$, and affine scaling gives
\begin{equation*}
 \norm{\nabla\tilde{v}}_{T_0\cup T_1}\le C\norm{\nabla v}_{T_0}.
\end{equation*}
Extending $\tilde{v}$ by zero on the remaining triangles gives a
continuous extension $\bar v$ to $\Omega_i$.

For the first assertion, write the interior patch cyclically as
$\Omega_i=\bigcup_{j=0}^{N-1}T_{j}$, where $T_{j}$ is bounded by the radial
edges $e_{j}$ and $e_{j+1}$, with indices modulo $N$. Subtract the mean of
$v$ on $T_0$ and denote the resulting polynomial again by $v$. Parameterize
the radial edges from $a_i$ to their outer endpoints by $t\in[0,1]$, and let
$g_0$ and $g_1$ be the traces of $v$ on $e_0$ and $e_1$.

For every $T_{j}$, the two affine vertex maps $T_{j}\to {T}_0$ produce copies
of $v$ whose ordered radial traces are $(g_0,g_1)$ and $(g_1,g_0)$.
Starting from $T_0$, choose these copies successively so that neighboring
traces agree. Each chosen copy $v_{j}$ satisfies
\begin{equation}\label{eq:patch-copy}
 v_{j}\in\Pp_r(T_{j}),\qquad
 \norm{\nabla v_{j}}_{T_{j}}\le C\norm{\nabla v}_{T_0}.
\end{equation}
If $N$ is even, this construction closes around the patch. If $N$ is odd,
apply it on $T_0,\ldots,T_{N-2}$. The two traces prescribed on the remaining
triangle $T_{N-1}$ are then the same function
$g\in\{g_0,g_1\}$. Their cross-edge seminorm satisfies
\begin{equation*}
 \int_0^1\int_0^1\frac{|g(s)-g(t)|^2}{(s+t)^2}\,ds\,dt
 \le |g|_{H^{1/2}(0,1)}^2.
\end{equation*}
On the reference triangle with radial edges $x=0$ and $y=0$,
the polynomial $g(x+y)$ has the prescribed traces.
Let $s\in[0,1]$ be the parameter measured from the common
vertex on each radial edge. Set $t=-s$ on one edge and $t=s$ on the other.
The two traces then form the single function $g(|t|)$ on $(-1,1)$.
The trace theorem on $T_0$, Poincar\'e's inequality, and
Corollary~\ref{cor:scalar-extension} therefore give a polynomial
$v_{N-1}\in\Pp_r(T_{N-1})$ with the prescribed traces and
\begin{equation*}
 \norm{\nabla v_{N-1}}_{T_{N-1}}
 \le C\norm{g(|\cdot|)}_{H^{1/2}(-1,1)}
 \le C\norm{g}_{H^{1/2}(0,1)}
 \le C\norm{\nabla v}_{T_0}.
\end{equation*}
Thus the pieces define a continuous extension and, by
\eqref{eq:patch-copy},
\begin{equation*}
 \norm{\nabla\bar v}_{\Omega_i}^2
 =\sum_{j=0}^{N-1}\norm{\nabla v_{j}}_{T_{j}}^2
 \le CN\norm{\nabla v}_{T_0}^2
 \le C\norm{\nabla v}_{T_0}^2.
\end{equation*}
The last inequality and all affine-scaling constants are uniform because
shape regularity bounds $N$. Restoring the subtracted mean completes the
proof.
\end{proof}

For test functions on the patch of degree at most $r$, with
$r\ge\ell+3$, use $\rho$ from
\eqref{eq:patch-ratio}.

\begin{proposition}\label{prop:patch}
Let $R$ be given by \eqref{eq:patch-res}, with $R(1)=0$ on an interior
patch. For $r\ge\ell+3$,
\begin{equation*}
 \dual R{X_i}\le C_{\rm patch}(\rho)\dual{R}{X_{i,r}},
\end{equation*}
where $C_{\rm patch}(\rho)$ depends only on shape regularity and an
upper bound for $\rho$ in \eqref{eq:patch-ratio}.
\end{proposition}

\begin{proof}
The reference estimates provide polynomial test functions
of degree at most $r$ for the volume corrections and edge transfers in the
vertex-patch reduction of
\cite[proof of Theorem~6.3]{CNSV17}.
We remove the volume residuals
first and then transfer the edge residuals around the patch.

\emph{Volume residuals.}
On each triangle, set $F_{T}=\phi_iG_T\in\Pp_\ell(T)$.
We construct a polynomial flux with divergence $F_T$ and
zero normal trace on the outer edge $o_T$. By
Theorem~\ref{thm:cm}, there is $\sigma_{T}^{(1)}\in\RT_\ell(T)$ satisfying
\begin{equation*}
 \diver\sigma_{T}^{(1)}=F_{T},\qquad
 \norm{\sigma_{T}^{(1)}}_{T}\le C\dual{F_{T}}{H^1_0(T)}.
\end{equation*}
Let $Y_{T}$ consist of $H^1(T)$ functions with zero trace on the two radial
edges. We define $w_{T}\in Y_{T}$ by
\begin{equation*}
 (\nabla w_{T},\nabla v)_{T}
 =-\langle\sigma_{T}^{(1)}\cdot n,v\rangle_{o_T},\qquad v\in Y_{T}.
\end{equation*}
Integration by parts gives
\begin{equation*}
 |\langle\sigma_{T}^{(1)}\cdot n,v\rangle_{o_T}|
 \le\bigl(\dual{F_{T}}{Y_{T}}+\norm{\sigma_{T}^{(1)}}_{T}\bigr)
 \norm{\nabla v}_{T}
 \le C\dual{F_{T}}{Y_{T}}\norm{\nabla v}_{T}.
\end{equation*}
It follows that $\tau_{T}=\nabla w_{T}$ is divergence-free, has normal
trace $-\sigma_{T}^{(1)}\cdot n$ on $o_T$, and satisfies
$\norm{\tau_{T}}_{T}\le C\dual{F_{T}}{Y_{T}}$.
This normal trace has degree at most $\ell$. Lemma~\ref{lem:solenoid}
therefore supplies a field $\sigma_{T}^{(2)}\in[\Pp_\ell(T)]^2$ with the
same normal trace and a comparable norm. Hence
$\sigma_{T}=\sigma_{T}^{(1)}+\sigma_{T}^{(2)}\in\RT_\ell(T)$ satisfies
\begin{equation}\label{eq:bulk-field}
 \diver\sigma_{T}=F_{T},\qquad \sigma_{T}\cdot n=0\text{ on }o_T,
 \qquad \norm{\sigma_{T}}_{T}\le C\dual{F_{T}}{Y_{T}}.
\end{equation}

For $v\in Y_T\cap\Pp_r(T)$, let $\bar v\in X_{i,r}$ be its
zero extension, with the patch mean subtracted on an interior patch.
By \eqref{eq:patch-res} and $R(1)=0$ on interior patches,
\begin{equation*}
 R(\bar v)=(F_T,v)_T,\qquad
 \norm{\nabla\bar v}_{\Omega_i}=\norm{\nabla v}_T.
\end{equation*}
Taking the supremum over these $v$ and applying Theorem~\ref{thm:mixed} gives
\begin{equation}\label{eq:bulk-delta}
  \dual{F_{T}}{Y_{T}}
 \le C\bigl(\mathfrak B_{\ell,r}+\mathfrak T_{\ell+1,r}\bigr)\dual{R}{X_{i,r}}
 \le C\rho^{3/4}\dual{R}{X_{i,r}}.
\end{equation}
We set $R_{T}(v)=(\sigma_{T},\nabla v)_{T}$. Integration by parts and
\eqref{eq:bulk-field} give
\begin{equation*}
 R_{T}(v)=-(F_{T},v)_{T}+
 \sum_{e\subset{\partial T\setminus o_T}}
 \langle\sigma_{T}\cdot n,v\rangle_e.
\end{equation*}
Thus $R_{T}$ cancels the volume functional and creates only radial-edge
functionals. Set $R^{(0)}=R+\sum_{T}R_{T}$.
Equations~\eqref{eq:bulk-field}--\eqref{eq:bulk-delta} give
\begin{equation}\label{eq:removed}
 \begin{aligned}
 \dual{R_{T}}{X_i}&\le C\rho^{3/4}\dual{R}{X_{i,r}},\\
 R^{(0)}(v)&=\sum_{e\text{ internal radial}}\int_e J_{i,e}^{(0)}v,
 \quad J_{i,e}^{(0)}\in\Pp_\ell(e).
 \end{aligned}
\end{equation}
On a boundary patch, densities on the boundary edges do not contribute
because the test functions vanish there. On an interior patch,
$R^{(0)}(1)=R(1)=0$. Since shape regularity bounds the number of
triangles, \eqref{eq:removed} implies
$\dual{R^{(0)}}{X_{i,r}}\le(1+C\rho^{3/4})\dual{R}{X_{i,r}}$.

\emph{Edge residuals.}
Label the triangles in an interior fan so that $T_{j}$ lies between radial
edges $e_{j}$ and $e_{j+1}$, with cyclic indices. Suppose that the current
density on $e_{j}$ is $J_{i,e_j}^{\rm cur}\in\Pp_\ell(e_{j})$. We solve the harmonic mixed
problem on $T_{j}$ with normal data $-J_{i,e_j}^{\rm cur}$ on $e_{j}$, zero normal data
on $o_{T_j}$, and zero Dirichlet data on $e_{j+1}$. Its gradient $\tau$
satisfies
\begin{equation*}
 \norm\tau_{T_{j}}=
 \left\|v\mapsto\int_{e_{j}}J_{i,e_j}^{\rm cur}v\right\|_{(H^1_{e_{j+1}}(T_{j}))'}.
\end{equation*}
Lemma~\ref{lem:solenoid}, applied to the two edges on which the normal trace is prescribed,
gives a divergence-free $\sigma_{j}\in[\Pp_\ell(T_{j})]^2$ with the same normal
data and a comparable norm. Adding $v\mapsto(\sigma_{j},\nabla v)_{T_{j}}$
cancels the density on $e_{j}$ and changes only that on $e_{j+1}$.
All densities retain degree at most $\ell$.

Theorem~\ref{thm:edge} gives polynomial tests of degree $r\ge\ell+2$ for the
mixed edge functional. By Lemma~\ref{lem:patch-extension}, each such
test extends to the adjacent triangle across $e_{j}$ and vanishes on the
two other radial edges of its support. Thus it pairs only with the
density being removed. 
At the $j$th step, let $R^{(j-1)}$ denote the current residual.
We obtain
\begin{equation*}
  \norm{\sigma_{j}}_{T_{j}}
 \le C\mathfrak T_{\ell,r}\dual{R^{(j-1)}}{X_{i,r}}
 \le C\rho^{3/4}\dual{R^{(j-1)}}{X_{i,r}}.
\end{equation*}
More precisely, let $\mathcal C_{j}(v)=(\sigma_{j},\nabla v)_{T_{j}}$ and
$R^{(j)}=R^{(j-1)}+\mathcal C_{j}$.
We have
\begin{equation}\label{eq:transfer-recurrence}
 \begin{aligned}
 \dual{\mathcal C_{j}}{X_i}&\le\norm{\sigma_{j}}_{T_{j}},\\
 \dual{R^{(j)}}{X_{i,r}}
 &\le \dual{R^{(j-1)}}{X_{i,r}}+\dual{\mathcal C_{j}}{X_{i,r}}
 &\le (1+C\rho^{3/4})\dual{R^{(j-1)}}{X_{i,r}}.
 \end{aligned}
\end{equation}
Mean subtraction on an interior patch leaves these pairings unchanged.

For a boundary fan with $N$ triangles,
use noncyclic indexing, with $e_1,e_{N+1}\subset\partial\Omega$.
Since test functions vanish on $e_1$, choose zero density there and set
$\sigma_1=0$, $\mathcal C_1=0$, and $R^{(1)}=R^{(0)}$.
Apply the transfer construction for $j=2,\ldots,N$.
Since the test functions vanish on $e_{N+1}$, we have
\begin{equation*}
 R^{(N)}(v)
 =\int_{e_{N+1}}(\sigma_N\cdot n_{T_N})\,v\,ds
 =0,\qquad \forall v\in X_i.
\end{equation*}

For an interior fan, we apply the transfers for $j=1,\ldots,N-2$.
The remaining densities lie on the two radial edges of $T_{N-1}$.
\mbox{Set $R^{\rm rem}=R^{(N-2)}$.}
Since $R^{(0)}(1)=0$ and $\mathcal C_{j}(1)=0$, we have
\begin{equation*}
 R^{\rm rem}(1)=R^{(0)}(1)+\sum_{j=1}^{N-2}\mathcal C_{j}(1)=0.
\end{equation*}
This value is the sum of the integrals of the two remaining edge densities.
We assign zero density to the outer edge of $T_{N-1}$ and apply
Theorem~\ref{thm:neumann}. Restriction of a patch test to $T_{N-1}$
bounds the remaining continuous patch norm by the triangle dual norm.
Conversely, Lemma~\ref{lem:patch-extension} extends every
polynomial test function on the triangle to the patch with a
controlled energy norm. If
$R^{\rm rem}$ is this remaining residual, these two observations give
\begin{equation*}
 \dual{R^{\rm rem}}{X_i}
 \le C\mathfrak T_{\ell,r}\dual{R^{\rm rem}}{X_{i,r}}
 \le C\rho^{3/4}\dual{R^{\rm rem}}{X_{i,r}}.
\end{equation*}

Set $R^{\rm rem}=R^{(N)}=0$ on a boundary patch. On either
type of patch, subtracting the volume and edge corrections gives
\begin{equation*}
 \dual R{X_i}\le
 \sum_T\dual{R_T}{X_i}+\sum_j\dual{\mathcal C_j}{X_i}
 +\dual{R^{\rm rem}}{X_i},
\end{equation*}
where the sums run over the patch triangles and the transfers performed.
By \eqref{eq:removed}, iteration of \eqref{eq:transfer-recurrence}, and
the bound on the final residual, the right-hand side is at most
$C_{\rm patch}(\rho)\dual R{X_{i,r}}$ since shape regularity bounds the
number of corrections.
The divergence and solenoidal liftings
have normal traces of degree at most $\ell$, so all reference
tests use the same local test degree $r\ge\ell+3$.
\end{proof}

\section{Global saturation and data oscillation}\label{sec:global}
\subsection{Polynomial-data saturation}
\begin{proof}[Proof of Theorem~\ref{thm:main}]
Fix $\lambda>1$, set $\ell=p$, and define
$\mathcal R_p(v)=(f,v)_\Omega-a(u_p,v)$ and its localization at a vertex $a_i$ by
$\mathcal{R}_i(v)=\mathcal R_p(\phi_iv)$ on $X_i$.

Elementwise integration by parts gives \eqref{eq:patch-res} with
\begin{equation*}
 \begin{aligned}
 G_T&=f+\Delta u_p\in\Pp_{\ell-1}(T),\\
 J_{i,e}&=-\phi_i\sum_{T:\,e\subset\partial {T}}
 \nabla(u_p|_{T})\cdot n_{T}\in\Pp_p(e)=\Pp_\ell(e),
 \end{aligned}
\end{equation*}
where $n_{T}$ is the outward unit normal of $T$.
On an interior patch,
$\mathcal{R}_i(1)=\mathcal R_p(\phi_i)=0$ since $\phi_i\in V_p$.

For $v\in H^1_0(\Omega)$, we take $c_{i}$ to be its patch mean at an
interior vertex and set $c_{i}=0$ at a boundary vertex. Partition of unity and Galerkin
orthogonality give
\begin{equation*}
 \mathcal R_p(v)=\sum_{i} \mathcal{R}_i(v-c_{i}).
\end{equation*}
By \eqref{eq:galerkin}, $\mathcal R_p(v)=a(u-u_p,v)$.
Each triangle belongs to three vertex patches, so
$\sum_i\norm{\nabla v}_{\Omega_i}^2=3\norm{\nabla v}_\Omega^2$.
Cauchy--Schwarz in the patch sum therefore gives
\begin{equation}\label{eq:localization}
 E_p=\dual{\mathcal R_p}{H^1_0(\Omega)}
 \le\sqrt3\left(\sum_{i}\dual{\mathcal{R}_i}{X_i}^2\right)^{1/2}.
\end{equation}
Let $q\ge\lceil\lambda p\rceil+3$ and
set $r=q-1$ as the polynomial degree bound for the test functions
on each patch, so multiplication by the piecewise affine hat function raises
degree at most to $q$. Then
\begin{equation*}
 r-\ell-2=q-p-3\ge(\lambda-1)p>0,\qquad
 \rho=1+\frac{p+3}{q-p-3}\le1+\frac4{\lambda-1}.
\end{equation*}
Since the degrees are integers, $r\ge\ell+3$. Hence
Proposition~\ref{prop:patch} has a constant depending only on $\lambda$
and shape regularity. If $v_{i}\in X_{i,r}$, zero extension gives
$\phi_iv_{i}\in V_q$. With $d=u_q-u_p$, it follows from \eqref{eq:galerkin} that
\begin{equation*}
 \mathcal{R}_i(v_{i})=a(d,\phi_iv_{i}).
\end{equation*}
Patch Poincar\'e's inequality, with zero mean or the boundary condition, and
$\norm{\nabla\phi_i}_{L^\infty}\le C/h_i$ imply
\begin{equation*}
 \norm{\nabla(\phi_iv_{i})}_{\Omega_i}\le C\norm{\nabla v_{i}}_{\Omega_i},\qquad
 \dual{\mathcal{R}_i}{X_{i,r}}\le C\norm{\nabla d}_{\Omega_i}.
\end{equation*}
We apply Proposition~\ref{prop:patch} in \eqref{eq:localization}
and use the discrete residual estimate above with
$\sum_i\norm{\nabla d}_{\Omega_i}^2=3\norm{\nabla d}_\Omega^2$. This gives
\begin{equation}\label{eq:increment}
 E_p\le \sqrt{3}\,C_{\rm patch}(\rho)
 \left(\sum_i\dual{\mathcal R_i}{X_{i,r}}^2\right)^{1/2}
 \le K_\lambda\norm{\nabla(u_q-u_p)}_\Omega.
\end{equation}
We may choose $K_\lambda>1$, independently of $p$, $q$, and the mesh sizes. Since $V_p\subset V_q$, Galerkin orthogonality gives
\begin{equation*}
 E_p^2=E_q^2+\norm{\nabla(u_q-u_p)}_\Omega^2.
\end{equation*}
Combining this identity with \eqref{eq:increment} gives
$E_q\le\alpha_\lambda E_p$ with $\alpha_\lambda=\sqrt{1-K_\lambda^{-2}}<1$.
This proves the theorem.
\end{proof}

\subsection{The oscillation estimate}
\label{sec:osc}
\begin{proof}[Proof of Theorem~\ref{thm:osc}]
Set $f_h|_T=\pi^{T}_{p-1}f$ and $\delta=\dual{f-f_h}{H^1_0(\Omega)}$.
The local approximation estimate \eqref{eq:local-projection} follows by
affine scaling from the reference-triangle case of \cite[Theorem~3.3]{BS00}.
For every $v\in H^1_0(\Omega)$, elementwise projection orthogonality,
weighted Cauchy--Schwarz, and \eqref{eq:local-projection} give
\begin{equation*}
 \begin{aligned}
 |(f-f_h,v)_\Omega|
 &=\left|\sum_{T\in\T}(f-f_h,v-\pi^T_{p-1}v)_T\right|\\
 &\le \osc_p(f,\T)
 \left(\sum_{T\in\T}\frac{p^2}{h_T^2}\norm{v-\pi^T_{p-1}v}_{T}^2\right)^{1/2}\\
 &\le C_{\rm osc}\osc_p(f,\T)\norm{\nabla v}_\Omega.
 \end{aligned}
\end{equation*}
Taking the supremum over $v$ of unit energy norm gives
$\delta\le C_{\rm osc}\osc_p(f,\T)$.

Let $\widetilde u$ and $\widetilde u_r$ solve the continuous and Galerkin
problems for $f_h$, and set
$\widetilde E_r=\norm{\nabla(\widetilde u-\widetilde u_r)}_\Omega$.
Subtracting the continuous equations in \eqref{eq:galerkin}
for $f$ and $f_h$ gives $a(u-\widetilde u,v)=(f-f_h,v)_\Omega$, so
$\norm{\nabla(u-\widetilde u)}_\Omega=\delta$.
If $\mathsf{P}^{\nabla}_r$ is the energy-orthogonal projector onto $V_r$,
Galerkin orthogonality gives $u_r=\mathsf{P}^{\nabla}_ru$ and
$\widetilde u_r=\mathsf{P}^{\nabla}_r\widetilde u$. Hence
\begin{equation*}
 (u-u_r)-(\widetilde u-\widetilde u_r)
 =(I-\mathsf{P}^{\nabla}_r)(u-\widetilde u).
\end{equation*}
Contractivity gives $|E_r-\widetilde E_r|\le\delta$.
Theorem~\ref{thm:main} applied to $f_h$ then implies
\begin{equation*}
 E_q\le\widetilde E_q+\delta
 \le\alpha_\lambda\widetilde E_p+\delta
 \le\alpha_\lambda E_p+(1+\alpha_\lambda)\delta.
\end{equation*}
Substituting the bound on $\delta$ gives the claimed estimate.
\end{proof}

	\section*{AI Disclosure}
	During the preparation of this work, GPT-5.6 sol and GPT-6.0 Astra (OpenAI) were used to explore proof strategies and assist with mathematical derivations. The authors reviewed and revised the AI-assisted material and independently verified all mathematical statements and proofs. The authors assume full responsibility for the results and the final manuscript.

\end{document}